\documentclass[12 pt,amsfonts] {amsart}
\font\azbuka=wncyr10
\def\cod{\mathop{\mathrm{cod}}\nolimits}
\def\dom{\mathop{\mathrm{dom}}\nolimits}
\def\id{\mathop{\mathrm{id}}\nolimits}
\def\R{{\Bbb R}}
\theoremstyle{definition}

\theoremstyle{remark}

\numberwithin{equation}{section}

\begin{document}

\title{ Flows with limited intersection of worldlines}

\author{Petra Augustov\'a (1)}
\address{(1) Institute of Information Theory and Automation of the Academy of Sciences of the Czech Republic, Prague}
\curraddr{}
\email{augustova@utia.cas.cz}
\thanks{The work of the first author is supported by Czech Science Foundation through the research grant no. 13-20433S}

\author{Lubom\'ir klapka}



\date{\today}

\dedicatory{}

\begin{abstract} In this paper we define a flow with
limited intersection of its worldlines and we construct and solve
functional equations for such flow using a special kind of set
embedding. For examples we use particular cases studied in the
past by different authors. The connection to higher order ordinary
differential equations is emphasized.
\end{abstract}

\maketitle

\section{Worldlines of the flow}
By the term {\it flow\/} we usually mean a one-parameter
transformation group (see e.g. Chapter 1 in \cite{01}, § 4 in
\cite{02} or section 1.3 v \cite{03}), i.e. the action $\R\times
M\ni(t,a)\mapsto ta\in M$ of the additive group of real numbers
$\R$ on some set $M$. The partial maps $\R\ni t\mapsto ta\in M$
defined for any $a\in\R$ are called {\it worldlines\/} of the
flow. It is clear that a flow is fully described by the set of its
worldlines.

The term worldline is connected with preference to space-time
approach in physics and is most used in relativity theory (e.g.
Sections 7 to 9 in \cite{04}). However its usage is not strictly
connected to this area only. In the theory of flows the so-called
{\it orbits\/} are intensively studied. Orbits are the sets
$\cup_{t\in\R}\{ta\}\subset M$, i.e. the images of worldlines.

Rigorously, the original Minkovski term worldline \cite{05} is not
the map $\R\to M$ directly but its graph. In the case of
worldlines $\R\ni t\mapsto ta\in M$ it means the set
$\cup_{t\in\R}\{(t,ta)\}\subset\R\times M$. This should not be
confusing since any map $\R\to M$ can be considered as a special
case of a binary relation on $\R\times M$ (see e.g. §3.4 in
\cite{06} or Section 6 of Chapter 1 in \cite{07}). Therefore it is
possible to identify the map with its graph and hence, $\R\ni
t\mapsto ta\in M$ and $\cup_{t\in\R}\{(t,ta)\}\subset\R\times M$
are just different descriptions of the same thing.

An important property of worldlines of one--parameter
transformation group is the fact that they are disjoint. This is
clear from the following reasoning. Two worldlines $\R\ni t\mapsto
ta_1\in M$ and $\R\ni t\mapsto ta_2\in M$ are different if and
only if there is $t_0\in\R$ where the two worldlines as maps do
not attain the same value, i.e. $t_0a_1\not=t_0a_2$. Since we can
cancel out the group element we successively obtain $a_1\not=a_2$
and $ta_1\not=ta_2$ for all $t\in\R$, so the worldlines are
disjoint as graphs.

For completeness let us note that the orbit invariance with
respect to the group action implies that the worldlines of the
flow $\R\times M\ni(t,a)\mapsto ta\in M$ are at the same time the
orbits of the flow  $\R\times(\R\times M)\ni(s,(t,a))\mapsto(s+t,
(s+t)a)\in(\R\times M)$.

\section{Generalization of flows}

We usually talk about flows in case of dynamical systems with
disjoint worldlines and each worldline can be uniquely determined
by anyone of its points. In this sense by flow we understand the
map that assigns to a point of a worldline the whole worldline.
Now we want to generalize the notion of the flow in case the
worldlines are not disjoint anymore but the number of
intersections of two different worldlines is further bounded above
by a natural number denoted by $k$. In this new situation we can
uniquely identify each worldline by its $k$ different points. Now
we can define the generalized flow as a map assigning to any
$k$-point restriction of any worldline the whole worldline. In the
sequel we formalize this construction.

Let $\R$ denote the set of all real numbers, $I\subset\R$ its
subset of cardinality $|I|>k$, $M$ any general given set, $M^I$
the set of all maps from $I$ to $M$. We define the {\it set of
worldlines with limited intersection} as the set  $X\subset M^I$
such that for any $x_1,x_2\in X$
\begin{equation}
\label{1} |x_1\cap x_2|\geq k\quad\Rightarrow\quad x_1=x_2.
\end{equation}
The cardinality of the intersection of maps $|x_1\cap x_2|$ is
defined here through identification of the maps with their graphs
as described above. From now on, by worldlines we understand the
elements of the set $X$. By the {\it set of all $k$-restrictions
of the worldline} we mean the set $C_k(x)$ of all
$k$-combinations, i.e. the set of all subsets $a\subset x$ such
that $|a|=k$. Since $x$ is a map then also any combination $a\in
C_k(x)$ is a map and satisfies $\dom a\in C_k(I)$, $a=x|_{\dom
a}$. We will call the set of all $k$-restrictions of all
worldlines $A=\cup_{x\in X}C_k(x)\subset\cup_{\alpha\in
C_k(I)}M^\alpha$ shortly as the {\it set of restrictions}.

We consider now the binary relation  $\varphi\subset A\times X$
defined by $(a,x)\in\varphi\,\Leftrightarrow\,a\subset x$. Since
$x\in X\subset M^I$ implies $|x|=|I|>k$, we have
$C_k(x)\not=\emptyset$. Hence, there exists $a\in A$ such that
$(a,x)\in\varphi$ and therefore $\varphi$ is surjective. By
definition of $A$, $\varphi$ is also co-surjective, i.e. for all
$a\in A$ there is $x\in X$ such that $(a,x)\in\varphi$. Moreover,
such $x$ is unique since $(a,x_1)\in\varphi\,
\wedge\,(a,x_2)\in\varphi$ yields $a\in C_k(x_1)\cap
C_k(x_2)=C_k(x_1\cap x_2)$, therefore $|x_1\cap x_2|\geq k$ and by
$(\ref{1})$ we get $x_1=x_2$. It follows that the relation
$\varphi$ is a map with domain $\dom\varphi=A$ and codomain
$\cod\varphi=X$.

This allows to define the {\it flow with limited intersection of
worldlines} as such a surjective map $\varphi\colon A\ni
a\mapsto\varphi_a \in X$, where $A=\cup_{x\in X}C_k(x)$ is the set
of restrictions, $X\subset M^I$ is the set of worldlines
satisfying the condition $(\ref{1})$, that $\varphi_a\supset a$
holds for all $a\in A$.

\section{Functional equations}

\noindent{\bf Theorem 1.} {\it The surjective map $\varphi\colon
a\mapsto\varphi_a$, $\dom\varphi\subset\cup_{\alpha\in
C_k(I)}M^\alpha$, $\cod \varphi\subset M^I$, is a flow with
limited intersection of worldlines if and only if for any
$a\in\dom\varphi$, $\beta\in C_k(I)$
\begin{equation}
\label{2}
\varphi_{\varphi_a|_\beta}=\varphi_a,\quad\quad\quad\varphi_a|_{
\dom a}=a.
\end{equation}}

\noindent{\bf Proof:} Let $\varphi$ be a flow with limited
intersection of worldlines and let $a\in\dom\varphi$, $\beta\in
C_k(I)$. Since $\dom\varphi_a=I$ we have $\varphi_a|_\beta\in
C_k(\varphi_a)\subset\dom\varphi$ and by the definition of flow
$\varphi_{ \varphi_a|_\beta}\supset\varphi_a|_\beta$. Because
$\varphi_a \supset\varphi_a|_\beta$ we obtain
$|\varphi_{\varphi_a|_\beta}
\cap\varphi_a|\geq|\varphi_a|_\beta|=|\beta|=k$ and using
$(\ref{1})$ we obtain $\varphi_{\varphi_a|_\beta}=\varphi_a$.
Furthermore, the definition of the flow guarantees also
$\varphi_a\supset a$ and hence, $\varphi_a|_{\dom a}\supset
a|_{\dom a}=a$. But $|\varphi_a|_{\dom a}|=|a|=k$ which implies
$\varphi_a|_{\dom a}=a$. Therefore, the equalities $(\ref{2})$ are
verified and the proof in one direction is finished.

Let now $\varphi$ be surjective with $\dom\varphi
\subset\cup_{\alpha\in C_k(I)}M^\alpha$ and $\cod\varphi\subset
M^I$ satisfying conditions $(\ref{2})$ for all $a\in\dom \varphi$
and $\beta\in C_k(I)$.  Then $\varphi_a|_{\dom a}=a$ implies
$\varphi_a\supset a$.

Let $x_1,x_2\in\cod\varphi$ satisfy $|x_1\cap x_2|\geq k$. By
surjectivity there exist $b,c\in\dom\varphi$ such that
$x_1=\varphi_b$, $x_2=\varphi_c$. The fact that intersection of
maps is a map again entails
$|\dom(\varphi_b\cap\varphi_c)|=|\varphi_b\cap\varphi_c|\geq k$,
therefore there exists $\delta\in
C_k(\dom(\varphi_b\cap\varphi_c))$. Then for such $\delta$ we have
$\varphi_b |_\delta=\varphi_c|_\delta$. By $(\ref{2})$, $x_1
=\varphi_b=\varphi_{\varphi_b|_\delta}=\varphi_{\varphi_c|
_\delta}=\varphi_c=x_2$. This proves condition $(\ref{1})$.

Assume that  $b\in\cup_{x\in\cod\varphi}C_k(x)$. Then there is
 $x_0\in\cod\varphi$ such that $b\in C_k(x_0)$. By surjectivity,
 there exist $a\in\dom\varphi$ such that $x_0=\varphi_a$.
We use the notation $\beta=\dom b$. Then, $\beta\in
C_k(\dom\varphi_a)=C_k(I)$. By $(\ref{2})$, we also have
$\varphi_{ \varphi_a|_\beta}=\varphi_a$, hence
$b=\varphi_a|_\beta\in\dom \varphi$. This shows
$\cup_{x\in\cod\varphi}C_k(x)\subset \dom\varphi$. Suppose on the
other hand that $a\in\dom\varphi$. Then $|a|=k$,
$a\subset\varphi_a$ and therefore $a\in C_k(\varphi_a)\subset
\cup_{x\in\cod\varphi}C_k(x)$. This proves $\dom\varphi=
\cup_{x\in\cod\varphi}C_k(x)$.

Because all the requirements of the definition are fulfilled, the
map $\varphi$ is a flow with limited intersection of worldlines.
This completes the proof in both directions. $\Box$

\section{Examples}

\noindent{\bf Flow without intersection of worldlines:} By
$(\ref{1})$, a flow without intersection of worldlines is a flow
with $k=1$. Usually, we assume the set of restrictions in the form
$\dom\varphi=\cup_{\alpha\in C_1(\R)}M^\alpha$. Among others this
also means that $a\in\dom\varphi$, $a=\{(r,m)\}$ where
$(r,m)\in\R\times M$ and $\cod\varphi\subset M^R$. We can rewrite
the functional equation $(\ref{2})$ for nonautonomous flow without
intersection of worldlines as
\begin{equation}
\label{3} F(t,r,F(r,s,m))=F(t,s,m),\quad F(s,s,m)=m,
\end{equation}
where $F\colon\R\times\R\times M\to M$ is a map with three
variables and the connection to the flow stands
$F(t,s,m)=\varphi_{\{(s,m)\}}(t)$.

The functional equations (\ref{3}) are the generalization to  
8.1.3(8)\cite{08} of the Sincov's equation (D. M. Sincov 1903). It
was already studied by many authors (A. N. Kolmogorov 1931, M.
Fr\'echet 1932, J. Acz\'el 1955, M. Hossz\'u 1958 and others) and
is used, i.a., in the theory of stochastic processes. Its detailed
description and references to the original work can be found in
8.1.3 and 8.1.4 of the Acz\'el's book \cite{08}.

\bigskip

\noindent{\bf Autonomous flow:} The most familiar flow without
intersection of worldlines is the autonomous flow.  The flow
autonomy means the closure property of $\cod\varphi$ under all
translation $\tau_s\colon\R\ni t\mapsto\, t+s \,\in\R$, i.e. $
x\in\cod\varphi\Rightarrow x\circ\tau_s\in\cod\varphi $ for all
$s\in\R$. By $(\ref{2})$, we know that
$\varphi_{\{(r,m)\}}(s)=\varphi_{\{(r,m)\}}\circ\tau_s(0)=$
$\varphi_{\varphi_{\{(r,m)\}}\circ\tau_s|_{\{r-s\}}}(0)=$
$\varphi_{\varphi_{\{(r,m)\}}|_{\{r\}}\circ\tau_s}(0)=$
$\varphi_{\varphi_{\{(r,m)\}}|_{\dom\{(r,m)\}}\circ\tau_s}(0)=$
$\varphi_{\{(r,m)\}\circ\tau_s}(0)=\varphi_{\{(r-s,m)\}}(0)$,
where $\circ$ means the composition of binary relations. Hence, we
can rewrite the functional equation $(\ref{2})$ of autonomous flow
without intersection of worldlines in the form
$$
F(r,F(s,m))=F(r+s,m),\quad F(0,m)=m,
$$
where $F\colon\R\times M\to M$ is map of two variables and its
relation to the flow are the following equalities
$F(r,m)=\varphi_{\{(r,m)\}}(0)$ and
$\varphi_{\{(r,m)\}}(s)=F(r-s,m)$.

The functional equations of autonomous flow without intersection
of worldlines describe a one-parametric group of transformations
$M\to M$ and are used under the name {\it translation equations\/}
in iteration theory and dynamical systems and in differentiable
case in the theory of Lie groups and first order differential
equations (see e.g. (3), (17) in \cite{09} and (1.9), (1.10) in
\cite{03}). The general solution is described in Section 8.2.2. of
\cite{08} and this book contains also references to the relevant
original works from the  middle of the last century.

\bigskip

\noindent{\bf The set of solutions of the Dirichlet problem in
ODE's:} We assume an ordinary differential equation $\ddot
x=f(t,x,\dot x)$ that has for any Dirichlet boundary condition
$x(\alpha)=a$, $x(\beta)=b$, $\alpha\not=\beta$ a unique solution
$x\colon\R\to\R^n$ satisfying this condition. The solutions of all
Dirichlet problems corresponding to such a differential equation
form a set of worldlines with limited intersections. Actually, the
worldlines that intersect at least twice satisfy the same
Dirichlet condition and hence they coincide. Thus, in this case
$k=2$, $I=\R$ and $M=\R^n$. This situation is described by P.
Chl\'adek in \cite{10} by using the equations
\begin{equation}
\label{4}
F(\tau,\alpha,\beta,a,b)=F(\tau,\gamma,\delta,F(\gamma,\alpha,\beta,a,b),
F(\delta,\alpha,\beta,a,b)),
\end{equation}
\begin{equation}
\label{5} F(\alpha,\alpha,\beta,a,b)=a,\quad
F(\beta,\alpha,\beta,a,b)=b,
\end{equation}
where $F\colon\R\times P_2(\R)\times(\R^n)^2\to\R^n$ is a function
and $P_2(\R)\subset\R^2$ the set of all $2$-permutations (without
repetition) of elements of $\R$. It is a special case of
functional equations $(\ref{2})$ where the unknown function $F$ is
related to the flow $\varphi\colon\cup_{\sigma\in
C_2(\R)}(\R^n)^\sigma\to\cod\varphi\subset(\R^n)^\R$ as follows:
\begin{equation}
\label{6}
F(\tau,\alpha,\beta,a,b)=\varphi_{\{(\alpha,a),(\beta,b)\}}(\tau).
\end{equation}

\bigskip

\noindent{\bf Harmonic oscillator:} A concrete case of the
previous example for $n=1$ and $k=2$ is the differential equation
$\ddot x=-x$ describing the behavior of a harmonic oscillator.
Since any of its maximal solutions has period $2\pi$ and any two
maximal solutions intersect, they must intersect infinitely many
times. Hence, they are not worldlines with limited intersection.
But we can consider the set $X$ as the set of all solutions
defined on an open interval $I$ of length at most $\pi$. Then the
condition $(\ref{1})$ is satisfied. To be exact, we fix
$I=(-\frac{\pi}{2},+\frac{\pi}{2})$. Then
$A=\dom\varphi=\cup_{\alpha\in C_2(I)}\R^\alpha$,
$X=\cod\varphi\subset\R^I$ and the equation $\varphi_a|_{\dom
a}=a$ admits for any $a\in\dom\varphi$ exactly one solution
$\varphi_a$ given by
$$
\varphi_a(t)=\sum_{(i,j)\in P_2(\dom
a)}a_i\,\,\,\frac{\sin(t-j)}{\sin(i-j)}.
$$
The relation $\varphi_{\varphi_a|_\beta}=\varphi_a$ becomes a
goniometrical identity where for any  $t\in I$, $a\in A$ and
$\beta\in C_2(I)$ we have
$$
\sum_{(r,s)\in P_2(\beta)}\sum_{(i,j)\in P_2(\dom
a)}a_i\frac{\sin(r-j)}{\sin(i-j)}\frac{\sin(t-s)}{\sin(r-s)}=\sum_{(i,j)\in
P_2(\dom a)}a_i\frac{\sin(t-j)}{\sin(i-j)}.
$$
Varying the interval $I$ the flow with limited intersection of
worldlines allows us to map step by step completely all the
solutions of the equation $\ddot x=-x$. In the last section of
this paper it is proven that similar process is possible for a
very wide class of ordinary differential equations.

\bigskip

\noindent{\bf Lagrange interpolating polynomials:} An interesting
example of a set of worldlines with limited intersection for a
general natural number $k$ is a set of polynomials. Let
$X=\cod\varphi\subset\R^\R$ is the set of all polynomials of
degree less than  $k$. The number $|x_1\cap x_2|$ for $x_1,x_2\in
X$ is the cardinality of the set of real roots of the polynomial
$x_1-x_2$. There are at least $k$ such roots if and only if we
deal here with the zero polynomial, i.e. $x_1=x_2$. The condition
$(\ref{1})$ is therefore satisfied. The set of restrictions is
$A=\dom\varphi=\cup_{\alpha\in C_k(\R)}\R^\alpha$ in this case.
The equation $\varphi_a|_{\dom a}=a$ has for any $a\in\dom\varphi$
exactly one solution $\varphi_a$, namely the {\it Lagrange
interpolation formula\/} (see e.\, g.  (2.5.3) in \cite{11}) given
by
$$
\varphi_a(t)=\sum_{i\in\dom a}a_i\prod_{j\in\dom
a\setminus\{i\}}\frac{t-j}{i-j}.
$$
The equation $\varphi_{\varphi_a|_\beta}=\varphi_a$ is then the
summation identity for the Lagrange interpolation formula, where
for any $t\in\R$, $a\in A$, $\beta\in C_k(\R)$ we have
$$
\sum_{r\in\beta}\sum_{i\in\dom a}a_i\prod_{j\in\dom
a\setminus\{i\}}\frac{r-j}{i-j}\prod_{s\in\beta\setminus\{r\}}\frac{t-s}{r-s}=\sum_{i\in\dom
a}a_i\prod_{j\in\dom a\setminus\{i\}}\frac{t-j}{i-j}.
$$
It is easy to see that the summation identity is equivalent to the
so called {\it Cauchy relations\/} (2.5.13) \cite{11}.

\section{Functional equation solutions}

\noindent In the sequel we use the term {\it $k$-frontal
embedding\/} of a set $W$ into the Cartesian power $M^q$ with
$q>k$ and by this we mean an embedding with the property that any
Cartesian projection of the set $W$ into each $M^k$ is a
bijection. Recalling the ordering of the set $I\subset\R$ we can
consider the sets $M^I$ and for $\alpha\subset I$ also $M^\alpha$
as Cartesian powers. Then by a $k$-frontal embedding into $M^I$ we
understand a surjection $\omega$ with $\cod\omega\subset M^I$ such
that any surjection
\begin{equation}
\label{7} \omega_\beta\colon\dom\omega\ni
w\mapsto\omega(w)|_\beta\in\cod \omega_\beta\subset M^\beta,
\end{equation}
where $\beta\in C_k(I)$, is a bijection at the same time. The
significance of this definition for the solution of functional
equations of a flow is clear from the following theorem.

\bigskip

\noindent{\bf Theorem 2.} {\it If $\omega$ is a $k$-frontal
embedding into $M^I$ then
\begin{equation}
\label{8} \varphi\colon\cup_{\sigma\in C_k(I)}\cod\omega_\sigma\ni
a\mapsto\omega\circ\omega_{\dom a}^{-1}(a)\in\cod\omega
\end{equation}
is a solution of functional equations $(\ref{2})$.

On the other hand, if the surjective map $\varphi$ with
$\dom\varphi\subset\cup_{\alpha\in C_k(I)}M^\alpha$ and $\cod
\varphi\subset M^I$ is a solution of functional equations
$(\ref{2})$ then there exists a $k$-frontal embedding into $M^I$
$\omega$ satisfying $(\ref{8})$.\/}

\bigskip

\noindent{\bf Proof:} Suppose that $\omega$ is a $k$-frontal
embedding into $M^I$ and $\varphi$ is a map given by $(\ref{8})$.
Then for any $a\in\dom\varphi$ and for any $\beta\in C_k(I)$ the
next equalities hold $\varphi_a|_\beta=(\omega\circ\omega_{\dom
a}^{-1}(a))|_\beta=\omega(\omega_{\dom
a}^{-1}(a))|_\beta=\omega_\beta\circ\omega_{\dom a}^{-1}(a)$.
Since $\dom\varphi_a|_\beta=\beta$ and thus $\varphi_a|_\beta\in
\cod\omega_\beta$, we get $\varphi_{\varphi_a|_\beta}=
\omega\circ\omega_\beta^{-1}({\omega_\beta\circ\omega_{\dom
a}^{-1}(a)})=\omega\circ\omega_{\dom a}^{-1}(a)=\varphi_a$,
$\varphi_a|_{\dom a}=\omega(\omega_{\dom a}^{-1}(a))|_{\dom
a}=\omega_{\dom a}(\omega_{\dom a}^{-1}(a))=a$, therefore
$\varphi$ is a solution of equations $(\ref{2})$.

Assume on the other side that $\varphi$ with
$\dom\varphi\subset\cup_{\alpha\in C_k(I)}M^\alpha$ and $\cod
\varphi\subset M^I$ is a solution of equations $(\ref{2})$. Let
$\omega=\id_{\cod\varphi}$. Then $\omega$ is a surjection and
$\cod\omega=\cod\varphi\subset M^I$. Consider  $\beta\in C_k(I)$
and two elements $w_1,w_2\in\dom\omega$ such that
$\omega_\beta(w_1)=\omega_\beta(w_2)$. Thanks to surjectivity of
the map $\varphi$ there exist  $a_1,a_2\in\dom\varphi$ such that
$\varphi_{a_1}=w_1$ and $\varphi_{a_2}=w_2$. Since
$\varphi_{a_1}|_\beta=\omega_\beta(w_1)=\omega_\beta(w_2)=\varphi_{a_2}|_\beta$,
by $(\ref{2})$ we also have
$w_1=\varphi_{a_1}=\varphi_{\varphi_{a_1}|_\beta}=\varphi_{\varphi_{a_2}|_\beta}=\varphi_{a_2}=w_2$.
This shows the injectivity of surjective maps $(\ref{7})$ and
therefore, their bijectivity. Hence,  $\omega$ is a $k$-frontal
embedding into $M^I$.

Since $\cod\omega_\sigma=\cup_{x\in\cod\varphi}\{x|_\sigma\}$ and
since $\cup_{\sigma\in C_k(I)}\{x|_\sigma\}=C_k(x)$, then
$\cup_{\sigma\in
C_k(I)}\cod\omega_\sigma=\cup_{x\in\cod\varphi}C_k(x)$ and by
Theorem 1 and by definition of a flow with limited intersection of
worldlines it follows that $\dom\varphi=\cup_{\sigma\in
C_k(I)}\cod\omega_\sigma$.

If we consider now $a\in\cup_{\sigma\in C_k(I)}\cod\omega_\sigma$,
by $(\ref{2})$, it holds $\cod\omega\ni\varphi_a=\omega_{\dom
a}^{-1}(\omega_{\dom a}(\varphi_a))=\omega_{\dom
a}^{-1}(\varphi_a|_{\dom a})=\omega_{\dom
a}^{-1}(a)=\omega\circ\omega_{\dom a}^{-1}(a)$, thus $\omega$
satisfies $(\ref{8})$. $\Box$

\bigskip

\noindent{\bf Solution for $k=1$:} We solve the equations
$(\ref{3})$. The embedding $\omega\colon M\to\cod\omega\subset
M^I$ can be expressed by the system  $\{G_t\}_{t\in I}$ of maps
$G_t\colon M\ni w\mapsto\omega(w)(t)\in M$. To obtain a
$1$-frontal embedding $\omega$ of $M$ into $M^I$ we need all the
maps $(\ref{7})$ to be bijections. This can be achieved only if
every map of the system $\{G_t\}_{t\in I}$ is a bijection. Indeed,
it is sufficient to take any system $\{G_t\}_{t\in I}$ of
bijections $G_t\colon M\to M$ and according to the Theorem 2 the map
$$F\colon\R\times\R\times M\ni(t,s,m)\mapsto
F(t,s,m)=\varphi_{\{(s,m)\}}(t)=$$$$\omega(\omega_{\{s\}}^{-1}(\{(s,m)\}))(t)=
G_t(\omega_{\{s\}}^{-1}(\{(s,m)\}))=G_t(G^{-1}_s(m))\in M$$ is a
solution of equations $(\ref{3})$. This solution is known already
for some time and is described by Theorem 2 in Section 8.1.3 of
the book \cite{08}.

\bigskip

\noindent{\bf Solution for $k=2$:} We can also solve the equations
$(\ref{4})$, $(\ref{5})$ thanks to a $2$-frontal embedding. We
express the embedding
$\omega\colon\R^{2n}\to\cod\omega\subset(\R^n)^\R$
 using maps of Cartesian spaces
$G\colon\R\times\R^{2n}\ni(\tau,w)\mapsto\omega(w)(\tau)\in\R^n$.
Should $\omega$ be a $2$-frontal embedding into $(\R^n)^\R$, the
equation
\begin{equation}
\label{9} G(\alpha,w)=a,\quad G(\beta,w)=b
\end{equation}
must have for any $(\alpha,\beta)\in P_2(\R)$,
$(a,b)\in\R^n\times\R^n$ exactly one solution for $w\in\R^{2n}$.
If we find a map $G$ with the properties above then at the same
time we have solution of equations $(\ref{4})$ and $(\ref{5})$.
Indeed, by $(\ref{6})$, $(\ref{7})$ and $(\ref{8})$ it holds
$F(\tau,\alpha,\beta,a,b)=G(\tau,H(\alpha,\beta,a,b))$ providing that
$H\colon P_2(\R)\times(\R^n)^2\to\R^{2n}$ is a map that assigns to
its arguments the solution
\begin{equation}
\label{10} w=H(\alpha,\beta,a,b)
\end{equation}
of equations $(\ref{9})$. In fact, elimitating $a$ and $b$ from
$(\ref{9})$ and $(\ref{10})$, we get
$w=H(\alpha,\beta,G(\alpha,w), G(\beta,w))$ which implies
$(\ref{4})$. On the other side, elimitating $w$ we get
$G(\alpha,H(\alpha,\beta,a,b))=a$,
$G(\beta,H(\alpha,\beta,a,b))=b$, which implies $(\ref{5})$. As
far as the authors of this paper are aware such general solution
of Chl\'adek's equations described herein are presented for the
first time.

\section{The smooth case}

For simplicity we restrict our attention to maps with continuous
partial derivatives of all orders and we call them {\it smooth\/}.
In the previous section we transferred the problem of solving
functional equations for a flow with limited intersection of
worldlines to the problem of finding the corresponding $k$-frontal
embedding. This arises the question of how difficult it is to find
such an embedding. This section shows that it is very easy in the
smooth case. The following Lemma gives sufficient condition for
local generating of $k$-frontal embedding by some smooth map.

\bigskip

\noindent{\bf Lemma.} {\it  Let $n$ and $k$ be natural numbers,
$G$ a smooth map, $\dom G\subset\R\times(\R^n)^k$ an open set,
$\cod G=\R^n$, $J\colon\dom G\to\R$ the Jacobian of the map
$$
\dom G\ni(t,w)\to\left(t,G(t,w),{\textstyle\frac{\partial
G(t,w)}{\partial
t}},\dots,{\textstyle\frac{\partial^{k-1}G(t,w)}{\partial
t^{k-1}}}\right)\in\R\times(\R^n)^k
$$
and let $(t_0,w_0)\in\dom G$ be a point where $J(t_0,w_0)\not=0$.
Then there exists such an open interval $I\subset\R$ and such an
open set  $U\subset(\R^n)^k$ that $(t_0,w_0)\in I\times
U\subset\dom G$ and the surjective map
\begin{equation}
\label{11} \omega\colon U\ni w\mapsto(I\ni t\mapsto
G(t,w)\in\R^n)\in\cod\omega
\end{equation}
is a $k$-frontal embedding of $U$ into $(\R^n)^I$.\/}

\bigskip

\noindent{\bf Proof:} Since $\dom G$ is an open set there is such
an open interval $I\subset\R$ and such an open set
$U\subset(\R^n)^k$ that $(t_0,w_0)\in I\times U\subset\dom G$. We
construct a map $K\colon I^k\times U\to(\R^n)^k$ assigning to any
$(\tau,w)\in I^k\times U$ the value $$
K(\tau,w)=(K_1(t_1,w),K_2(t_1,t_2,w),\dots,K_k(t_1,t_2,\dots,t_k,w)),
$$
where $\tau=(t_1,t_2,\dots,t_k)$, $K_1(t_1,w)=G(t_1,w)$ and $K_2,
\dots, K_k$ are recurrently given by
$$K_{i+1}(t_1,t_2,\dots,t_{i+1},w)=
\int_0^1 K'_i(t_1(1-s)+t_2\,s,\,t_3,\dots,t_{i+1},w) \,\,\,ds,$$
where $K'_i$ denotes the partial derivative of $K_i$ with respect
to the first variable. A direct calculation of the integral on the
right side easily verifies that the smoothness of $K_i$ implies
the smoothness of $K_{i+1}$. The map $K_1=G|_{I\times U}$ is
smooth by the assumptions of the Lemma, hence the map $K$ as
constructed above is also smooth. In the point $(\tau_0,w_0)$,
where $\tau_0=(t_0,t_0,\dots,t_0)\in I^k$, the Jacobian of the
smooth map $(\tau,w)\mapsto(\tau,K(\tau,w))$ is equal to $J(t_0,w_0)$
and therefore, by assumptions of the Lemma is nonzero.

Using the inverse function theorem (see e.g. Theorem 1A.1 in
\cite{12}) we can shrink the open interval $I$ an the open set $U$
preserving the condition $(t_0,w_0)\in I\times U\subset\dom G$ in
such a way that $(\tau,w)\mapsto(\tau,K(\tau,w))$ is an injection.
But then also the map $w\mapsto K(\tau,w)$ is an injection for any
$\tau\in I^k$. We consider $\beta\in C_k(I)$ and a k-permutation
without repetition $(t_1,t_2,\dots,t_k)\in
P_k(\beta)\subset\beta^k\subset I^k$. Since, for $i=1,2,\dots,k$
we have
$$
K_i(t_1,t_2,\dots,t_i,w)=\sum_{j\in\{1,2,\dots,i\}}G(t_j,w)
\prod_{m\in\{1,2,\dots,i\}\setminus\{j\}}\frac{1} {t_j-t_m},
$$
the map $w\mapsto(G(t_1,w),G(t_2,w),\dots,G(t_k,w))$ is an
injection too. Finally, by $(\ref{7})$ and $(\ref{11})$, for any
$\beta\in C_k(I)$ the surjection
$$
\omega_\beta\colon
w\mapsto\{(t_1,G(t_1,w)),(t_2,G(t_2,w)),\dots,(t_k,G(t_k,w))\}
$$
is also an injection, i.e. $\omega_\beta$ is a bijection. By
definition of $k$-frontal embedding the map $(\ref{11})$ is a
$k$-frontal embedding of $U$ into $(\R^n)^I$ and the proof is
complete. $\Box$

\section{Ordinary differential equations}

Now we formulate and prove a theorem on local description of the
set of solutions of a differential equation using a flow with
limited intersection of worldlines. Consider an ordinary
differential equation
\begin{equation}
\label{12}
x^{(k)}(t)=f(t,x^{(0)}(t),x^{(1)}(t),\dots,x^{(k-1)}(t)),
\end{equation}
where $t$ is an independent variable,  $x^{(i)}$ are derivatives
of the unknown function $x$, especially $x^{(0)}=x$ and
$x^{(i+1)}(t)=dx^{(i)}(t)/dt$, with Cauchy conditions
\begin{equation}
\label{13} (x^{(0)}(t_0),x^{(1)}(t_0),\dots,x^{(k-1)}(t_0))=w.
\end{equation}

\bigskip

\noindent{\bf Theorem 3.} {\it Let $n$ and $k$ be natural numbers,
let $f$ be a smooth map where $\dom f\subset\R\times(\R^n)^k$ is
an open set and $\cod f=\R^n$ and let $(t_0,w_0)\in\dom f$. Then
there exists an open interval $I\subset\R$ and an open set
$U\subset(\R^n)^k$ such that $(t_0,w_0)\in I\times U\subset\dom f$
and the set $X$ of all solutions $x$ of  the ordinary differential
equation $(\ref{12})$ defined on $I$ and satisfying the Cauchy
conditions $(\ref{13})$ with $w\in U$ is a set of worldlines with
limited intersection.}

\bigskip

\noindent{\bf Proof:} We make use of the Lemma. For $\dom
G\subset\R\times(\R^n)^k$ we consider the set of all pairs $(t,w)$
such that $(t_0,w)\in\dom f$ and $t\in\dom x$ with $x$ being the
maximal solution of the ordinary differential equation
$(\ref{12})$ satisfying the Cauchy conditions $(\ref{13})$ where
$w\in U$. For $\cod G$ take $\R^n$ and set
\begin{equation}
\label{14} G(t,w)=x(t).
\end{equation}
Clearly,  $(t_0,w_0)\in\dom G$. The variable $y\colon\dom x\ni
t\mapsto(x^{(0)}(t),$ $x^{(1)}(t),\dots,x^{(k-1)}(t))\in(\R^n)^k$
satisfies the first order ordinary differential equation
$dy(t)/dt=\Psi(t,y(t))$, where $\Psi$ is the smooth map
$(t,(a_0,a_1,$
$\dots,a_{k-1}))\mapsto(a_1,a_2,\dots,a_{k-1},f(t,(a_0,a_1,\dots,
a_{k-1})))$ defined on the open set $\dom\Psi=\dom f$. Every
maximal solution of this differential equation satisfying the
initial condition $y(t_0)=y_0$ is in the form $y\colon
t\mapsto\Phi(t,t_0,y_0)$, where $\Phi$ is a continuous map defined
on an open set $\dom\Phi\subset\R\times\dom f$ (see Theorem 14,
§23 in \cite{13}). Since the map $\Phi$ is also smooth (see
Corollary 4, §7 in \cite{02}) and since $G(t,w)=p(\Phi(t,t_0,
w))$, where $p\colon(\R^n)^k\ni(a_0,a_1,\dots,a_{k-1})\mapsto
a_0\in\R^n$ denotes the Cartesian projection on the zero factor,
then also $G$ is a smooth map defined on an open set $\dom
G=\{(t,w)|(t,t_0,w)\in\dom\Phi\}\subset\R\times(\R^n)^k$. If we
calculate the value of the Jacobian following the guidelines from
the Lemma we get $J(t_0,w_0)=1$, hence all the requirements of the
Lemma are fulfilled.

In this case the conclusions of the Lemma are accomplished too,
i.e. there exist such an open interval $I\subset\R$ and such open
sets $U\subset(\R^n)^k$ that $(t_0,w_0)\in I\times U\subset\dom
G\subset\dom f$ and there is a $k$-frontal embedding $\omega$
given by $(\ref{11})$. By $(\ref{14})$ the set $\cod\omega$ is the
set of all solutions  $x$ of the ordinary differential equation
$(\ref{12})$ defined on the interval $I$ and satisfying the Cauchy
conditions $(\ref{13})$ with $w\in U$.

By Theorem 2, the map $\varphi$ defined in $(\ref{8})$ is a
solution of functional equations $(\ref{2})$. By Theorem 1,
$\varphi$ is indeed a flow with limited intersection of
worldlines. The definition of flow with limited intersection of
worldlines says that the set $X=\cod\varphi$ is the set of
worldlines satisfying the condition $(\ref{1})$. The definition of
the set of worldlines with limited intersection then identifies
the set $X$ as such a set. Finally, $(\ref{8})$ gives that
$X=\cod\varphi=\cod\omega$ and this fact finishes the proof.
$\Box$


\begin{thebibliography}{99}
\footnotesize
\bibitem{01}
Auslander, J., {\it Minimal flows and their extensions\/},
North-Holland Mathematical Studies, vol. 153, North-Holland,
Amsterdam, New York, Oxford and Tokyo, 1988.
\bibitem{02}
Arnold, V. I., {\it Ordinary differential equations}, New York:
Springer-Verlag, 1992.
\bibitem{03}
Olver, P., {\it Applications of Lie groups to differential
equations}, (2nd edition), New York: Springer-Verlag, 1993.
\bibitem{04}
Felice, F. De, Clarke, C. J., {\it Relativity on Curved
Manifolds\/}, Cambridge University Press, 1992.
\bibitem{05}
Minkowski, H., {\it Raum und Zeit\/}, aus: Jahresberichte der
Deutschen Mathematiker-Vereinigung, Leipzig, 1909, Teubner,
Leipzig 1909.
\bibitem{06}
Suppes, P., {\it Axiomatic set theory\/}, Van Nostrand, Princeton,
1960.
\bibitem{07}
Levy, A., {\it Basic Set Theory\/}, Springer-Verlag, Berlin, 1979.
\bibitem{08}
Acz\'el, J., {\it Lectures on functional equations and their
applications}, Academic press, New York and London, 1966.
\bibitem{09}
Balibrea, F., Reich, L. and Smital, J., {\it Iteration Theory:
Dynamical Systems and Functional Equations\/}, Int. J. Bifurcation
Chaos 13 (2003), 1627--1647.
\bibitem{10}
Chl\'adek, P.,{\it The functional formulation of second-order
ordinary differential equations}, Aeq. Math. 69 (2005), No. 3,
263--270.
\bibitem{11}
Davis, P. J., {\it Interpolation and Approximation\/}, Dover, New
York, 1975.
\bibitem{12}
Dontchev, A. L., and Rockafellar, R. T., {\it Implicit Functions
and Solution Mappings\/}, Springer, New York, 2009.
\bibitem{13}
{\azbuka Pontryagin, L. S., {\emph Obyknovennye differencialnye
uravneniya}, Moskva, Nauka, 1974.}

\end{thebibliography}
\end{document}